\documentclass{article}

\usepackage{amsmath}

\newtheorem{theorem}{Theorem} 

\newtheorem{ theorem}{Theorem}

\newtheorem{lemma}{Lemma} 

\newtheorem{ lemma}{Lemma}

\newtheorem{fact}{Fact}

\newtheorem{ conjecture}{Conjecture}

\newtheorem{ proposition}{Proposition}

\title{{\Large \textbf{On Pure Degree Sequence Manipulations \\ Forcing Long Cycles in Graphs}}}
\author{\normalsize  Zhora Nikoghosyan}

\begin{document}
\count0=1

\maketitle

\begin{abstract}

The well-known Hamiltonian sufficient conditions, proposed by Dirac, Faudree et al., P\'osa, Bondy, Chv\'atal are based on pure degree manipulations without any additional conditions. In this paper, we present two new types of pure degree manipulations that produce relatively simpler (in view of verification) results. The reverse versions (long-cycle versions) of the obtained results are presented as well. The importance of pure degree manipulations  is motivated by their amenability (as starting points) to generalizations and new ideas in a great  variety of ways.  

Keywords: Hamilton cycle, longest cycle, circumference, minimum degree, degree sequence.

\end {abstract}

\section{Introduction}

We consider only finite undirected graphs with neither loops nor multiple edges. A good reference for any undefined terms is \cite{[4]}.

The set of vertices of a graph $G$ is denoted by $V(G)$ and the set of edges - by $E(G)$. Let $n$ be the order (the number of vertices) of $G$, $c$ the order of a longest cycle and $p$ the order of a longest path in $G$. For a given vertex $v$, we use $N(v)$ to denote the set of all neighbors of $v$ and $d(v)=|N(v)|$ to denote the degree of $v$. Let $d_1,d_2,...,d_n$ be the degree sequence of $G$ with $d_1\le d_2\le...\le d_n$, where the first term $d_1$ is known as the minimum degree of $G$ denoted by $\delta$.  

A graph $G$ is Hamiltonian if $G$ contains a Hamilton cycle, that is a simple spanning cycle. We will call $n$-closure of $G$ the graph obtained by recursively joining pairs $x,y$ of nonadjacent vertices such that $d(x)+d(y)\geq n$ until no such pair remains. It will be denoted by $Cl_n(G)$.

We write a cycle (a path) $Q$ with a given orientation by $\overrightarrow{Q}$. For $x\in V(Q)$, we denote the successor and the predecessor of $x$ on $\overrightarrow{Q}$ (if such vertices exist) by $x^+$ and $x^-$, respectively. For $U\subseteq V(Q)$, we denote $U^+=\{u^+|u\in U\}$ and $U^-=\{u^-|u\in U\}$. We use $P=x\overrightarrow{P}y$ to denote a path with end vertices $x$ and $y$ in the direction from $x$ to $y$. We say that a vertex $z_1$ precedes vertex $z_2$ on $\overrightarrow{Q}$ if $z_1$, $z_2$ occur on $\overrightarrow{Q}$ in that order, and indicate this relationship by $z_1\prec z_2$. We will write $z_1\preceq z_2$ when either $z_1=z_2$ or $z_1\prec z_2$.

In 1952, Dirac \cite{[6]} proved that a graph $G$ of order $n$ is Hamiltonian if every vertex of $G$ has degree at least $\frac{n}{2}$. 
The bound $\delta\ge\frac{n}{2}$ in this theorem is sharp in the sense that it cannot be weakened to $\delta\ge\frac{n-1}{2}$. On the other hand, the same theorem alternatively can be formulated as follows:

\begin{theorem} (Dirac \cite{[6]}, 1952). If $d_1\ge \frac{n}{2}$ in a graph $G$ with degree sequence $d_1\le d_2\le \ldots \le d_n$, then $G$  is Hamiltonian.
\end{theorem}

Formulated in this way, Dirac's theorem is far from to be best possible and it is evident in which direction it can be extended. 

Historically, along with vertex degrees, a number of other graph characteristics have been appeared due to their impact on graph structures, such as size (the number of edges), edge and vertex connectivity, binding number, toughness, independence number, forbidden subgraphs, and so on. 

However, some well-known extensions of Dirac's theorem, proposed by Faudree et al., P\'osa, Bondy and Chv\'atal are purely based on degree manipulations without any additional conditions.

Following this direction, in 2007, Faudree \cite{[7]} proved that every graph of order $n$ is Hamiltonian if $|\{v\in V(G)|d(v)<\frac{n}{2}\}|\leq\delta-1.$ It is not hard to see that this result can be equivalently reformulated \cite{[8]} in terms of a much simpler algebraic relation.

\begin{theorem}  (Faudree et al. \cite{[7]},\cite{[8]}, 2007). Let $d_1\le d_2\le \ldots \le d_n$ be the degree sequence of a graph $G$. If $d_\delta\ge \frac{n}{2}$, then $G$ is Hamiltonian.
\end{theorem}

The next extension of Dirac's theorem through pure degree manipulations was given by P\'osa \cite{[11]}.

\begin{theorem} (P\'osa \cite{[11]}, 1962). Let $d_1\le d_2\le \ldots \le d_n$ be the degree sequence of a graph $G$. If $d_i\ge i+1$ for each $i<\frac{n}{2}$, then $G$ is Hamiltonian.
\end{theorem}

In 1971,  Bondy \cite{[1]} obtained a new result based on degree sums.

\begin{theorem} (Bondy \cite{[1]}, 1971). Let $d_1\le d_2\le \ldots \le d_n$ be the degree sequence of a graph $G$. If $d_i+d_j\ge n$ whenever $i<j$, $d_i\le i$, $d_j\le j-1$ for each $i, j$, then $G$ is  Hamiltonian.
\end{theorem}

Finally, Chv\'atal \cite{[5]} obtained a P\'osa-type result by manipulating the degrees in the segment $\frac{n}{2}\leq i\leq n$: if $d_{n-i}\ge n-i$ for each positive integer $i<\frac{n}{2}$, then $G$ is Hamiltonian. In fact, by combining this result with P\'osa's theorem, Chv\'atal proved a more general result.

\begin{theorem} (Chv\'atal \cite{[5]}, 1972).
Let $d_1\le d_2\le \ldots \le d_n$ be the degree sequence of a graph $G$. If either $d_i\ge i+1$ or $d_{n-i}\ge n-i$ for each positive integer $i<\frac{n}{2}$, then $G$ is Hamiltonian.
\end{theorem}

In this paper we present two new pure degree manipulations forcing the graph to be Hamiltonian.

\begin{ theorem}. Let $d_1\le d_2\le \ldots \le d_n$ be the degree sequence of a graph $G$. If $d_\delta+d_{\delta+1}\ge n$, then $G$ is Hamiltonian. 
\end{ theorem}

\begin{ theorem}. Let $d_1\le d_2\le \ldots \le d_n$ be the degree sequence of a graph $G$. If $\min\{d_\delta+d_{\delta+2},2d_{\delta+1}\}\ge n$, then $G$ is Hamiltonian. 
\end{ theorem}

In the next section, we will show that Theorem 1 and Theorem 2 are strictly stronger than Dirac's Theorem A and are independent of Ore's and P\'osa's theorems.

Among the non-pure degree manipulations, the most striking example is the famous Ore's theorem \cite{[10]} - a direct extension of Dirac's theorem. Additional limitations in Ore's theorem based on the degree sums of nonadjacent vertices significantly change the nature of the impact on graph structures to another level.  

\begin{theorem} (Ore \cite{[10]}, 1960). If $d(u)+d(v)\ge n$ for each pair of nonadjacent vertices $u,v$ of a graph $G$ of order $n$, then $G$ is Hamiltonian.
\end{theorem}

The next Ore-type condition was given by Faudree \cite{[7]} in terms of the number of all missing edges $xy$ with $d(x)+d(y)<n$. Let

$$
g(n,\delta)=\begin{cases}
\infty, & \text{if \ $n\leq 2\delta$};\\
\frac{n^2-1}{8}, & \text{if \ $2\delta+1\leq n\leq 6\delta-3$ and $n$ is odd};\\
\frac{n^2+2n-8}{8}, & \text{if \ $2\delta+2\leq n\leq 6\delta-4$ and $n$ is even};\\
\delta n-\frac{3}{2}\delta^2-\frac{1}{2}\delta, & \text{if $n\geq 6\delta-2$}. 
\end{cases}
$$

Further, let
$$
N_2(G)=\{(x,y)|xy\notin E(G)  \ \text{and}  \ d(x)+d(y)<n \}; \ n_2(G)=|N_2(G)|.
$$

\begin{theorem} (Faudree et al. \cite{[7]}, 2007).
Every graph of order $n$ is Hamiltonian if $n_2(G)<g(n,\delta)$.
\end{theorem}

In Section 2, we will show that Theorem G is stronger than Ore's theorem but not stronger than P\'osa's theorem. 

A new Ore-type condition was proposed by Las Vergnas \cite{[9]}.

\begin{theorem} (Las Vergnas \cite{[9]}, 1971).
Let $G$ be a graph with vertices $v_1$, $v_2$,..., $v_n$ and let there be no $i,j$ with
\begin{itemize}
\item $i<j, \  v_iv_j\notin E(G),$ 
\item $d(v_i)\le i,\  d(v_j)\le j-1,$
\item $d(v_i)+d(v_j)\le n-1, \  i+j\ge n.$
\end{itemize}
Then $G$ is Hamiltonian.
\end{theorem}

Finally, Bondy and Chv\'atal \cite{[2]}, using the degree sums and nonadjacent vertices with a recursively vertices joining procedure, obtained the following power result inspired by Ore's theorem.

\begin{theorem} (Bondy-Chv\'atal \cite{[2]}, 1976)
A graph of order $n$ is Hamiltonian if and only if its $n$-closure is Hamiltonian.
\end{theorem}

In the next section we will present the relationship between all above theorems.

The reverse versions (the long-cycle versions) of Theorem B, Theorem 1 and Theorem 2 can be developed around the reverse version of Dirac's next theorem \cite{[6]} on long cycles for 2-connected graphs: in every 2-connected graph, $c\ge \min\{p,2\delta\}$. This theorem in terms of degree sequences is equivalent to the following.

\begin{theorem} (Dirac \cite{[6]}, 1952). Let $G$ be a 2-connected graph with degree sequence $d_1\le d_2\le \ldots \le d_n$. Then $c\ge \min\{p,2d_1\}$. 
\end{theorem}

The following extension of Theorem J was given in \cite{[8]}.

\begin{theorem} (\cite{[8]}, 2017). Let $G$ be a 2-connected graph with degree sequence $d_1\le d_2\le \ldots \le d_n$.  Then $c\ge \min\{p,2d_\delta\}$. 
\end{theorem}

In this paper we present the following two further extensions of Theorem K.

\begin{ theorem}. Let $d_1\le d_2\le \ldots \le d_n$ be the degree sequence of a 2-connected graph $G$. Then $c\ge \min\{p, d_\delta+d_{\delta+1}\}$. 
\end{ theorem}

\begin{ theorem}. Let $d_1\le d_2\le \ldots \le d_n$ be the degree sequence of a 2-connected graph $G$. Then $c\ge \min\{p,  d_\delta+d_{\delta+2}, 2d_{\delta+1}\}$. \end{ theorem}

To show that Theorem 1 is sharp, let $G=K_\delta+(\delta K_1\cup K_2)$, where $d_{\delta+1}=\delta+1$ and $n=2\delta+2$. Since $G$ is a non-Hamiltonian graph with $2d_{\delta+1}\ge n$, the condition $d_\delta+d_{\delta+1}\ge n$ in Theorem 1 cannot be relaxed to $d_{\delta+1}+d_{\delta+1}\ge n$.

Next, if $G=K_{\delta}+\overline{K}_{\delta+1}$, then $d_\delta=d_{\delta+1}=\delta$, $d_{\delta+2}=2\delta$ and $n=2\delta+1$. Since $G$ is a non-Hamiltonian graph with $d_\delta+d_{\delta+2}\ge n$ and $d_\delta+d_{\delta+1}\ge n-1$,  the condition $d_\delta+d_{\delta+1}\ge n$ in Theorem 1 cannot be relaxed to $d_\delta+d_{\delta+2}\ge n$ or $d_\delta+d_{\delta+1}\ge n-1$.

Thus Theorem 1 is best possible in all respects.

To show that Theorem 2 is also sharp in all respects, suppose first that $G=K_{\delta}+\overline{K}_{\delta+1}$, where $d_\delta=d_{\delta+1}=\delta$, $d_{\delta+2}=2\delta$ and $n=2\delta+1$. Since $G$ is a non-Hamiltonian graph with
\begin{align*}  
(a1)  \qquad & \min\{d_\delta+d_{\delta+2},2d_{\delta+1}+1\}\ge n, \qquad \qquad \qquad \qquad \qquad\qquad \qquad \\
(a2) \qquad &  \min\{d_\delta+d_{\delta+2},2d_{\delta+2}\}\ge n, \qquad \qquad \qquad \qquad \qquad \qquad \qquad
\end{align*}
the condition $\min\{d_\delta+d_{\delta+2},2d_{\delta+1}\}\ge n$ in Theorem 2 cannot be relaxed to (a1) or (a2).

If $G=K_\delta+(\delta K_1\cup K_2)$, then $d_\delta=\delta$, $d_{\delta+1}=d_{\delta+2}=\delta+1$, $d_{\delta+3}=2\delta+1$ and $n=2\delta+2$. Since
$G$ is a non-Hamiltonian graph with
\begin{align*}
(a3) \qquad & \min\{d_\delta+d_{\delta+2},2d_{\delta+1}\}\ge n, \qquad \qquad \qquad \qquad \qquad\qquad \qquad \\
(a4)  \qquad  &  \min\{d_{\delta+1}+d_{\delta+2},2d_{\delta+1}\}\ge n, \qquad \qquad \qquad \qquad \qquad\qquad \qquad  \\
(a5)  \qquad  &  \min\{d_\delta+d_{\delta+3},2d_{\delta+1}\}\ge n, \qquad \qquad \qquad \qquad \qquad\qquad \qquad   
\end{align*}
the condition $\min\{d_\delta+d_{\delta+2},2d_{\delta+1}\}\ge n$ in Theorem 2 cannot be relaxed to (a3), (a4) or (a5).

The bound $c\ge \min\{p,d_\delta+d_{\delta+1}\}$ in Theorem 3 cannot be strengthened to $c\ge \min\{p+1,d_\delta+d_{\delta+1}\}$, since the condition $c\le p$ holds in any case.

Let $\kappa$ be the vertex connectivity of $G$. If $G=K_\delta+(\delta K_1\cup K_2)$, then 
$$
\kappa\ge 2, \ d_{\delta+1}=\delta+1, \ c=2\delta+1, \  p=2\delta+2.
$$

Since $c< \min\{p,d_{\delta+1}+d_{\delta+1}\}$, the bound $c\ge \min\{p,d_\delta+d_{\delta+1}\}$ in Theorem 3 cannot be strengthened to $c\ge \min\{p,d_{\delta+1}+d_{\delta+1}\}$.

If $G=K_{\delta}+\overline{K}_{\delta+1}$, then 
$$
\kappa\ge 2,\   d_\delta=d_{\delta+1}=\delta, \ d_{\delta+2}=2\delta, \ c=2\delta, \  p=2\delta+1.
$$

Since 
$$c<\min\{p,d_\delta+d_{\delta+2}\} \ \  \text{and} \ \  c<\min\{p,d_\delta+d_{\delta+1}+1\},
$$
the bound $c\ge \min\{p,d_\delta+d_{\delta+1}\}$ in Theorem 3 cannot be strengthened to $c\ge \min\{p,d_\delta+d_{\delta+2}\}$ or $c\ge \min\{p,d_\delta+d_{\delta+1}+1\}$.

Finally, if $G=K_1+2K_\delta$, then $\kappa=1$, $d_\delta=d_{\delta+1}=\delta$, $c=\delta+1$ and $p=2\delta+1$. Since $c<\min\{p,d_\delta+d_{\delta+1}\}$, the 2-connectivity condition in Theorem 3 cannot be relaxed.

Thus Theorem 3 is best possible in all respects.

Now we turn to the sharpness of Theorem 4. The bound $c\ge p$ cannot be strengthened to $c\ge p+1$, since the condition $c\le p$ holds in any case.

If $G=K_{\delta}+\overline{K}_{\delta+1}$, then 
$$
\kappa\ge 2, \ d_\delta=d_{\delta+1}=\delta, \  d_{\delta+2}=2\delta,  \ c=2\delta, \ p=2\delta+1.
$$

Since
\begin{align}
c&<\min\{p,d_\delta+d_{\delta+2},2d_{\delta+1}+1\}, \nonumber \\
c&<\min\{p,d_\delta+d_{\delta+2},2d_{\delta+2}\}, \nonumber
\end{align}
the bound $c\ge\min\{p,d_\delta+d_{\delta+2},2d_{\delta+1}\}$ in Theorem 4 cannot be strengthened to $c\ge\min\{p,d_\delta+d_{\delta+2},2d_{\delta+1}+1\}$ or $c\ge\min\{p,d_\delta+d_{\delta+2},2d_{\delta+2}\}$.

If $G=K_\delta+(\delta K_1\cup K_2)$, then 
$$
\kappa\ge 2, \  d_\delta=\delta, \ 
d_{\delta+1}=d_{\delta+2}=\delta+1,
$$
$$
d_{\delta+3}=2\delta+1, \ c=2\delta+1, \ p=2\delta+2.
$$

Since
\begin{align}
c&<\min\{p,d_{\delta+1}+d_{\delta+2},2d_{\delta+1}\}, \nonumber \\
c&<\min\{p,d_\delta+d_{\delta+3},2d_{\delta+1}\}, \nonumber
\end{align}
the bound $c\ge\min\{p,d_\delta+d_{\delta+2},2d_{\delta+1}\}$ in Theorem 4 cannot be strengthened to $c\ge\min\{p,d_{\delta+1}+d_{\delta+2},2d_{\delta+1}\}$ or $c\ge\min\{p,d_\delta+d_{\delta+3},2d_{\delta+1}\}$.

Finally if $G=K_1+2K_\delta$, then 
$$
d_\delta=d_{\delta+1}=d_{\delta+2}=\delta, \  c=\delta+1, \ p=2\delta+1.
$$

Since $\kappa=1$ and $c< \min\{p,d_\delta+d_{\delta+2},2d_{\delta+1}\}$, the 2-connectivity condition in Theorem 4 cannot be relaxed.

Thus Theorem 4 is sharp in all respects.

\section{Relationships between theorems}

Recall the following fact which includes a number of well-known Hamiltonian degree conditions.

\begin{fact}.
Each second theorem among \{Theorem A, Theorem F, Theorem C, Theorem D, Theorem E, Theorem H, Theorem I\} is strictly stronger than the previous one.
\end{fact}

Let $C_i$ denote the cycle $v_1v_2...v_iv_1$ for each $i\in\{1,2,...,n\}$.

Theorem $G$ vs Ore's theorem and P\'osa's theorem: Let $\pi_0=(3^2,4,5^4,6^3)$ be the degree sequence of a graph $H_{0}$ obtained from $C_{10}$ by adding the edges $v_1v_3$, $v_2v_4$, $v_2v_5$, $v_2v_8$, $v_3v_5$, $v_3v_6$, $v_3v_8$, $v_4v_6$, $v_4v_7$, $v_4v_8$, $v_5v_7$, $v_6v_9$ and $v_7v_9$. We have 
$$
n=10, \ \delta=3 \ \  \text{and} \ \  2\delta+2\leq n \leq 6\delta-4.
$$

In addition, we have 
$$
n_2(H_0)=14 \ \ \text{and} \ \ g(n,\delta)=\frac{1}{8}(n^2+2n-8)=14.
$$

Since $n_2(G)=g(n,\delta)$, the sequence $\pi_0$ fails to satisfy $n_2(G)<g(n,\delta)$. On the other hand, $\pi_0$ explicitly satisfies the condition of P\'osa's theorem. So we have the following.

\begin{fact}.
Theorem G is stronger than Ore's theorem but not stronger than P\'osa's theorem.
\end{fact}

Let us adjust the relationships between Theorem A, Theorem B, Theorem 1 and Theorem 2. Notice first that 
$$
\min\{d_{\delta}+d_{\delta+2}, 2d_{d_{\delta+1}}\}\geq d_\delta+d_{\delta+1}\geq 2d_\delta\geq 2d_1=2\delta.
$$

It immediately follows from these observations that Theorem 2 is stronger than Theorem 1; Theorem 1 is stronger than Theorem B; and Theorem B is stronger than Dirac's theorem.

Now let's construct some examples of graphs to show that there are strictly strong connections between these four theorems. 

Theorem A vs Theorem B: Let $H_1$ be a graph with degree sequence $\pi_1=(2,3,3,3,3)$, or briefly $\pi_1=(2,3^4)$, obtained from $C_5$ by adding the edges $v_2v_4$ and $v_3v_5$. Then  
$$
d_\delta =d_2=3>\frac{5}{2}=\frac{n}{2} \ \  \text{and} \ \ \delta=2<\frac{5}{2}=\frac{n}{2},
$$
implying that $\pi_1$ satisfies the condition of Theorem B and fails to satisfy the condition of Dirac's theorem. In other words, Theorem B is strictly stronger than Dirac's theorem. 

Theorem 1 vs Theorem B: Denote by $H_2$ the graph with degree sequence $\pi_2=(2,3,4^2,5^3)$ obtained from $C_7$ by adding the edges $v_2v_4$, $v_2v_5$, $v_2v_6$, $v_3v_5$, $v_3v_6$, $v_3v_7$ and $v_4v_6$. Then 
$$
d_\delta+d_{\delta+1}=7=n \ \ \text{and} \ \ d_\delta=d_2=3<\frac{7}{2}=\frac{n}{2}
$$
i.e., $\pi_2$ satisfies the condition of Theorem 1 and fails to satisfy the condition of Theorem B. This means that Theorem 1 is strictly stronger than Theorem B.

Theorem 2 vs Theorem 1: Let $H_3$ be a graph with degree sequence $\pi_3=(2,3,4,5^5)$ obtained from $C_8$ by adding the edges $v_2v_5$, $v_3v_6$, $v_3v_7$, $v_3v_8$, $v_4v_6$, $v_4v_7$, $v_4v_8$, $v_5v_7$ and $v_5v_8$. For  $\pi_3$, we have 
$$
\min\{d_\delta+d_{\delta+2},2d_{\delta+1}\}=8=n \ \ \text{and} \ \  d_\delta+d_{\delta+1}=7<n.
$$

Therefore $\pi_3$ satisfies the condition of Theorem 2 and fails to satisfy the condition of Theorem 1. It follows that Theorem 2 is strictly stronger than Theorem 1. Thus we have the following.

\begin{fact}.
Theorem 2 is strictly stronger than Theorem 1; Theorem 1 is strictly stronger than Theorem B; and Theorem B is strictly stronger than Theorem A.
\end{fact}

Now let's see how Ore's theorem is related to Theorem B, Theorem 1, and Theorem 2. 

Theorem B vs Ore's theorem: Take disjoint $K_4$ and $K_7$  with vertex sets 
$$
V(K_4)=\{x_1,x_2,x_3,x_4\}, \ V(K_7)=\{y_1,y_2,...,y_7\}.
$$

Let $H_4$ be a graph on 11 vertices with degree sequence $\pi_4=(4,5^3,7^7)$ obtained from $K_4$ and $K_7$ by adding the edges $x_iy_i$ $(i=1,2,3,4)$, $x_2y_5$, $x_3y_6$ and $x_4y_7.$
Clearly, 
$$
d(x_1)=4,\ d(x_i)=5 \ \ (i=2,3,4), \ \  d(y_i)=7 \ \ (i=1,2,...,7),
$$
$$
\text{and} \ \ d_\delta=d_4=5<\frac{11}{2}=\frac{n}{2}.
$$

It follows that $\pi_4$ satisfies the condition of Ore's theorem and fails to satisfy the condition of Theorem B.

To show the opposite, let $H_5$ be a graph with degree sequence $\pi_5=(2,3^2,4^3)$ obtained from $C_6$ by adding the edges $v_2v_4$, $v_2v_5$, $v_3v_6$ and $v_4v_6$. Clearly, 
$$
d_\delta=d_2=3=\frac{n}{2} \ \  \text{and} \ \  d(v_1)+d(v_3)=5<6=n,
$$
implying that $\pi_5$ satisfies the condition of Theorem B and fails to satisfy the condition of Ore's theorem.

Thus Theorem B and Ore's Theorem are mutually independent.

Theorem 1 vs Theorem 2: Let $K_9$ and $K_{16}$ be disjoint complete graphs with vertex sets $V(K_9)=\{x_0,x_1,...,x_8\}$ and $V(K_{16})=\{y_1,y_2,...,y_{16}\}$. Further, let $H_6$ be a graph of order 25 obtained from $K_9$ and $K_{16}$ by adding the following edges
 $$
x_iy_i \ (i=1,2,...,8), \ x_1y_8, \ x_{i+1}y_i
\ (i=1,2,...,7), 
$$
$$
x_iy_{i+8} \  (i=1,2,...,8), \ x_1y_{16}, \ \ x_iy_{i+7} \  (i=2,3,...,8).
$$

Since $\pi_6=(8,12^8,17^{16})$
is the degree sequence of $H_6$, it is easy to see that $d(u)+d(v)=8+17=25=n$ for each nonadjacent vertices $u,v\in V(H_6)$. Next we have 
$$
d_\delta+d_{\delta+1}=24<n \ \  \text{and} \ \ \min\{d_\delta+d_{\delta+2}, 2d_{\delta+1}\}=24<n.
$$

In other words,  $\pi_6$ satisfies the condition of Ore's theorem and fails to satisfy the conditions of Theorem B, Theorem 1 and Theorem 2. 

To show the opposite, let $H_7$ be a graph with degree sequence $\pi_7=(3^3,5^5)$ obtained from $C_8$ by adding the edges $v_2v_4$, $v_3v_1$, $v_3v_7$, $v_3v_8$, $v_4v_1$, $v_4v_6$, $v_5v_8$, $v_6v_8$ and $v_6v_1$. Clearly, 
$$
d_\delta+d_{\delta+1}=8=n, \ \ \min\{d_\delta+d_{\delta+2},
2d_{\delta+1}\}=8=n 
$$
$$
\text{and} \ \  d(v_2)+d(v_5)=6<8=n.
$$

This means that  $\pi_7$ satisfies the conditions of Theorem 1 and Theorem 2; and fails to satisfy the condition of Ore's theorem.

Thus Theorem 1 and Theorem 2 are independent of Ore's theorem. Recalling the above observations, we can claim the following.

\begin{fact}.
Theorem B, Theorem 1 and Theorem 2 are independent of Ore's theorem.
\end{fact}

Now let's see how P\'osa's theorem relates to Theorem B, Theorem 1 and Theorem 2. 

P\'osa's theorem vs Theorem B:  To prove that P\'osa's theorem is stronger than Theorem B, let's show that if $d_\delta\geq \frac{n}{2}$, then $d_i\geq i+1$ for each $i<\frac{n}{2}$. Assume the contrary, that is $d_i\leq i$ for some $i<\frac{n}{2}$. If $\delta\leq i$, then $d_\delta\leq d_i \leq i<\frac{n}{2}$, a contradiction. Hence $\delta \geq i+1$. On the other hand, from $d_i\leq i$ we have $\delta=d_1\leq d_i\leq i$, which yields $\delta\leq i$. This contradiction shows that P\'osa's theorem is indeed stronger than Theorem B.

Furthermore, let $H_8$ be a graph with degree sequence $\pi_8=(2,3,4^2,5^3)$ obtained from $C_7$ by adding the edges $v_2v_4$, $v_2v_5$, $v_2v_6$, $v_3v_5$, $v_3v_6$, $v_3v_7$ and $v_4v_6$. Then 
$$
d_i\geq i+1 \ \  (i=1,2,3) \ \ \text{and} \ \  d_\delta=d_2=3<\frac{7}{2}=\frac{n}{2}.
$$

This means that $\pi_8$ satisfies the condition of P\'osa's theorem and fails to satisfy the condition $d_\delta\geq \frac{n}{2}$. Therefore P\'osa's theorem is strictly stronger than Theorem B.

P\'osa's theorem vs Theorem 1: Let $H_9$ be a graph with degree sequence $\pi_9=(3^3,5^5)$ obtained from $C_8$ by adding the edges $v_1v_3$, $v_2v_6$, $v_2v_7$, $v_2v_8$, $v_3v_5$, $v_3v_8$, $v_4v_7$, $v_5v_7$ and $v_5v_8$. We have $d_\delta+d_{\delta+1}=8=n$ and $d_\delta= \delta=3<\frac{n}{2}$, that is $\pi_9$ satisfies the condition of Theorem 1 and fails to satisfy the condition of P\'osa's theorem. 

To show the opposite, let $\pi_{10}=(2,3,4^5,5)$ be the  degree sequence of a graph $H_{10}$ obtained from $C_8$ by adding the edges $v_2v_5$, $v_2v_7$, $v_3v_5$, $v_3v_6$, $v_4v_6$, $v_4v_7$ and $v_6v_8$. Then $d_i\geq i+1$ for $i=1,2,3$ and $d_\delta+d_{\delta+1}=7<8$, that is $\pi_{10}$ satisfies the condition of P\'osa's theorem and fails to satisfy the condition $d_\delta+d_{\delta+1}\geq n$. Hence Theorem 1 and P\'osa's theorem are independent.

Theorem 2 vs P\'osa's theorem: Let $H_{11}$ be a graph with degree sequence $\pi_{11}=(3^3,5^5)$ obtained from $C_8$ by adding the edges $v_1v_3$, $v_2v_6$, $v_2v_7$, $v_2v_8$, $v_3v_5$, $v_3v_8$, $v_4v_7$, $v_5v_7$ and $v_5v_8$. We have 
$$
\min\{d_\delta+d_{\delta+2},2d_{\delta+1}\}=8=n \ \ \text{and} \ \  d_\delta=d_3=3,
$$
which means that $\pi_{11}$ satisfies the condition of Theorem 2 and fails to satisfy the condition of P\'osa's theorem. 

To show the opposite, let $\pi_{0}=(3^2,4,5^4,6^3)$ be the degree sequence of  $H_{12}$. Then 
$$
d_i\geq i+1 \ \  (i=1,2,3,4) \ \ \text{and} \ \ \min\{d_\delta+d_{\delta+2},2d_{\delta+1}\}=9<n=10,
$$
which means that $\pi_{0}$ satisfies the condition of P\'osa's theorem and  fails to satisfy the condition of Theorem 2. Thus Theorem 2 and theorem of P\'osa are independent as well.

Combining above observations, we can claim the following.

\begin{fact}.
Theorem 1 and Theorem 2 are independent of P\'osa's theorem; P\'osa's theorem is strictly stronger than Theorem B.
\end{fact}

To complete this section, we shall show that Bondy's theorem is strictly stronger than Theorem 2, implying in particular that Bondy's theorem is strictly stronger than Theorem 1 and Theorem B, as well. The fact that Bondy's theorem stronger than Theorem 2, is shown in the proofs of Theorem 1 and Theorem 2. 

For strictly superiority, take $H_{0}$  with degree sequence $\pi_{0}=(3^2,4,5^4,6^3)$. Then the minimum $i$ with $d_i\leq i$ is $i=5$. Further the minimum $j$ with $d_j\leq j-1$ is $j=6$. Therefore $d_i+d_j=10= n$ for each $i, j$ with $d_i\leq i$ and $d_j\leq j-1$. This means that $\pi_{0}$ satisfies the condition of Bondy's theorem. On the other hand, $\pi_{0}$ fails to satisfy the condition of Theorem 2 since 
$$
\min\{d_\delta+d_{\delta+2},2d_{\delta+1}\}=9<10=n.
$$

Thus Bondy's theorem is strictly stronger than Theorem 2.

Due to above observations, we have the following.

\begin{fact}.
Bondy's theorem is strictly stronger than Theorem B, Theorem 1 and Theorem 2.
\end{fact}

\section{Preliminaries}

Let $\overrightarrow{P}=v_1v_2...v_p$ be a longest path in $G$. Clearly, $N(v_1)\cup N(v_p)\subseteq V(P)$. A vine of length $m$ on $P$ is a set
$$
\{L_i=w_i\overrightarrow{L}_iz_i: 1\le i\le m\}
$$
of internally-disjoint paths such that
\begin{itemize}
    \item $V(L_i)\cap V(P)=\{w_i,z_i\} \ \ (i=1,...,m)$,
    \item $v_1=w_1\prec w_2\prec z_1\preceq w_3\prec z_2\preceq w_4\prec ...\preceq w_m\prec z_{m-1}\prec z_m=v_p$ \mbox{on} $P$.
\end{itemize}

The following result, known as the vine Lemma, guarantees the existence of at least one vine on $\overrightarrow{P}$ in a 2-connected graph (the original vine's idea was developed in \cite{[6]}).

\begin{lemma} (Bondy, Locke \cite{[3]}, 1981). Let $G$ be a $k$-connected graph and $P$ a path in $G$. Then there are $k-1$ pairwise-disjoint vines on $P$.
\end{lemma}

\section{Proofs}

\textbf{Proof of Theorem 1}. Let $d_1\le d_2\le \ldots \le d_n$ be the degree sequence of a graph $G$ and let $d_\delta+d_{\delta+1}\geq n$. To show that $G$ is Hamiltonian, by Theorem D (Bondy), it suffices to show that $d_i+d_j\geq n$ whenever $d_i\leq i$ and $d_j\leq j-1$.   
If  $i<\delta$, then $i\geq d_i\leq i\leq \delta-1=d_1-1$, a contradiction. Let $i\geq \delta$, implying that $d_i\geq d_\delta$. Next, if $j\leq \delta$, then $d_j\leq j-1\leq \delta-1=d_1-1$, again a contradiction. Let $j\geq \delta+1$, implying that $d_j\geq d_{\delta+1}$.   Then 
$$
d_i+d_j\geq d_{\delta}+d_{\delta+1}\geq n
$$
and we are done by Theorem D. 
Theorem 1 is proved. \qquad \rule{7pt}{6pt}

\textbf{Proof of Theorem 2}. Let $d_1\le d_2\le \ldots \le d_n$ be the degree sequence of a graph $G$ and let $\min\{d_\delta+d_{\delta+2}, 2d_{\delta+1}\}\geq n$. By Theorem D, it suffices to show that $d_i+d_j\geq n$ whenever $d_i\leq i$ and $d_j\leq j-1$. 

As in the proof of Theorem 1, we can show that  $d_i\geq d_\delta$ and $d_j\geq d_{\delta+1}$.  Further, if $i\geq\delta+1$, that is $d_i\geq d_{\delta+1}$, then 
$$
d_i+d_j\geq 2d_{\delta+1}\geq\min\{d_\delta+d_{\delta+2}, 2d_{\delta+1}\}\geq n
$$
and we are done. 
Let $i=\delta$. Since $d_i\leq i$, i.e., $d_\delta \leq \delta$, we have $d_\delta=\delta$. If $j\geq \delta+2$, then 
$$
d_i+d_j\geq d_{\delta}+d_{\delta+2}\geq\min\{d_\delta+d_{\delta+2}, 2d_{\delta+1}\}\geq n
$$
and we are done.
Let $j=\delta+1$. Similarly, since $d_j\leq j-1$, we have $d_{\delta+1}\leq \delta$, which yields $d_{\delta+1}=\delta$. So, 
$i=\delta$, $j=\delta+1$ and $d_\delta=d_{\delta+1}=\delta$, implying that
$$
d_i+d_j= d_{\delta}+d_{\delta+1}=2d_{\delta+1}\geq \min\{d_\delta+d_{\delta+2}, 2d_{\delta+1}\}\geq n.
$$

By Theorem D, $G$ is Hamiltonian. \qquad \rule{7pt}{6pt}

\textbf{Proof of Theorem 3}. Let $G$ be a 2-connected graph and let $\overrightarrow{P}=v_1v_2...v_p$ be a longest path in $G$. Clearly, $N(v_1)\cup N(v_p)\subseteq V(P)$. Assume that
\begin{align*}
 (a1) \quad  \ P \ \text{is chosen so that} \ d(v_1) \ \text{is maximum. \qquad \qquad \qquad \qquad \qquad \qquad}  
\end{align*}

Let $x_1,x_2,...,x_t$ be the elements of $N(v_1)$ occurring on $\overrightarrow{P}$ in a consecutive order, where $t=d(v_1)\ge \delta$. Observe that for each $i\in \{1,2,...,t\}$,
$$
x_i^-\overleftarrow{P}v_1x_i\overrightarrow{P}v_p
$$
is a longest path in $G$, implying that
$$
N(x_i^-)\subseteq V(P) \quad  (i=1,2,...,t).
$$
By (a1),
\begin{equation}\label{1}
 d(v_1)\ge d(x_i^-) \quad (i=1,2,..,t), \quad d(v_1)\ge d(v_p).
\end{equation}

If $x_t=v_p$, then $c\ge p+1$. Let $x_t\not=v_p$, that is $x_t\prec v_p$.

Further assume that
\begin{align*}
 (a2) \quad  \ P \ \text{is chosen so that} \ d(v_p) \ \text{is maximum, subject to (a1)}. \qquad \qquad   
 \end{align*}
 
Let $y_1,y_2,...,y_f$ be the elements of $N(v_p)$ occurring on $\overleftarrow{P}$ in a consecutive order. Since
$$
y_i^+\overrightarrow{P}v_py_i\overleftarrow{P}v_1
$$
is a longest path in $G$, we have
$$
N(y_i^+)\subseteq V(P) \quad  (i=1,2,...,f).
$$
By (a2),
\begin{equation}\label{2}
 d(v_p)\ge d(y_i^+) \quad (i=1,2,..,f).
\end{equation}

Further, by (1) and (2),

\begin{equation*}
    \begin{split}
    d(v_1)&\ge \max\{d(x_1^-),d(x_2^-),...,d(x_t^-),d(v_p)\}\\
&\ge\max\{d_1,d_2,...,d_{t+1}\}=d_{t+1}= d_{d(v_1)+1}
    \end{split}
\end{equation*}
and
\begin{equation*}
\begin{split}
d(v_p)&\ge \max\{d(y_1^+),d(y_2^+),...,d(y_f^+\}\\
&\ge \max\{d_1,d_2,...,d_f\}=d_f= d_{d(v_p)},
\end{split}
\end{equation*}
implying that

\begin{equation}\label{3}
d(v_1)+d(v_p)\ge d_{d(v_1)+1}+d_{d(v_p)}\ge d_\delta+d_{\delta+1}.
\end{equation}

\textbf{Case 1}. $x_t\preceq y_f$.

Let
$$
\{L_i=w_i\overrightarrow{L}_iz_i: 1\le i\le m\}
$$
 be a vine of minimal length $m$ on $\overrightarrow{P}$. Since $P$ is a longest path in $G$, we have $L_1,L_M\in E(G)$. Next since $m$ is minimal, we have $x_t\prec z_2$, $x_t\prec w_3$ and $w_{m-1}\prec y_f$, $z_{m-2}\prec y_f$. Choose $z_1^*\in V(P)$ such that $w_2\prec z_1^*$ and $|V(w_2\overrightarrow{P}z_1^*)|$ is minimal. Analogously, choose $w_m^*\in V(P)$ such that $w_m^*\prec z_{m-1}$ and $|V(w_m^*\overrightarrow{P}z_{m-1})|$ is minimal. Put
$$
H=P\cup \bigcup_{i=2}^{m-1}L_i\cup \{v_1z_1^*, v_pw_m^*\}.
$$
By deleting the following paths
$$
w_i\overrightarrow{P}z_{i-1} \quad (i=3,4,...,m-1), \quad w_2\overrightarrow{P}z_1^*, \quad w_m^*\overrightarrow{P}z_{m-1}
$$
from $H$ (except for their end vertices), we obtain a cycle $C$ with at least $d(v_1)+d(v_p)+1$ vertices. By (3),
$$
c\ge |V(C)|\ge d(v_1)+d(v_p)+1\ge d_\delta+d_{\delta+1}.
$$

\textbf{Case 2}. $y_f\prec x_t$.

\textbf{Case 2.1}. $N(v_1)\cap N^+(v_p)\not=\emptyset$.

Let $v\in N(v_1)\cap N^+(v_p)$, that is $v_1v,v_pv^-\in E(G)$. Clearly,
$$
v_1v\overrightarrow{P}v_pv^-\overleftarrow{P}v_1
$$
is a cycle of order $p$, that is $c\ge p$.

\textbf{Case 2.2}. $N(v_1)\cap N^+(v_p)=\emptyset$.

Since $y_f\prec x_t$, we can choose $x_i\in N(v_1)$ and $y_j\in N(v_p)$ such that $y_j\prec x_i$ and $v_1v,v_pv\not\in E(G)$ for each vertex $v$ with $y_j\prec v\prec x_i$. Put
$$
C=v_1x_i\overrightarrow{P}v_py_i\overleftarrow{P}v_1.
$$

Then
\begin{equation*}
\begin{split}
c&\ge |V(C)|\ge |N(v_1)|+|N^+(v_p)|+|\{v_1\}|-|\{y_j\}|\\
&\ge|N(v_1)|+|N(v_p)|=d(v_1)+d(v_p).
\end{split}    
\end{equation*}

By (3), $c\ge d_\delta+d_{\delta+1}$.  \qquad \rule{7pt}{6pt}

\noindent\textbf{Proof of Theorem 4}. Let $G$ be a 2-connected graph and let $\overrightarrow{P}=v_1v_2...v_p$ be a longest path in $G$. Define the vertices $x_1,x_2,...,x_t,y_1,y_2,...,y_f$ as in proof of Theorem 1, where we have proved

\begin{equation}
d(v_1)\ge d_{d(v_1)+1}, \quad  d(v_p)\ge d_{d(v_p)}.
\end{equation}

\textbf{Case 1}. $x_t\preceq y_f$.

As in proof of Theorem 1 (Case 1), we can form a cycle with at least $d(v_1)+d(v_p)+1$ vertices.  By (4),
\begin{equation}
c\ge d(v_1)+d(v_p)+1\ge d_{d(v_1)+1}+d_{d(v_p)}+1.
\end{equation}
If $d(v_p)\ge \delta+1$, then by (5), $c\ge 2d_{\delta+1}$ and we are done. Otherwise we have $d(v_p)= \delta$, implying that
$$
d(y^+_1)=d(y^+_2)=...=d(y^+_f)=\delta.
$$

If $d(v_1)=\delta$, then clearly, $d(v_p)=d(v_1)$ and
\begin{equation*}
    \begin{split}
     d(v_p)&\ge \max\{d(y_1^+),d(y_2^+),...,d(y_f^+),d(v_1)\}\\  
    & \ge \max\{d_1,d_2,...,d_{f+1},\}=d_{f+1}=d_{d(v_p)+1}\ge d_{\delta+1}.
    \end{split}
\end{equation*}

Next by (4), $d(v_1)\ge d_{d(v_1)+1}\ge d_{\delta+1}$. Thus by (5),
$$
c\ge d(v_1)+d(v_p)+1\ge 2d_{\delta+1}+1.
$$
Let $d(v_1)\ge \delta+1$. By (4), $d(v_1)\ge d_{d(v_1)+1}\ge d_{\delta+2}$ and $d(v_p)\ge d_{d(v_p)}\ge d_\delta$. Further, by (5),
$$
c\ge d(v_1)+d(v_p)+1\ge d_{\delta+2}+d_\delta+1,
$$
and again we are done.

\textbf{Case 2}. $y_f\prec x_t$.

\textbf{Case 2.1}. $N(v_1)\cap N^+(v_p)\not=\emptyset$.

Let $v\in N(v_1)\cap N^+(v_p)$, that is $v_1v,v_pv^-\in E(G)$. Clearly,
$$
v_1v\overrightarrow{P}v_pv^-\overleftarrow{P}v_1
$$
is a cycle of order $p$, that is $c\ge p$. 

\textbf{Case 2.2}. $N(v_1)\cap N^+(v_p)=\emptyset$.

Since $y_f\prec x_t$, we can choose $x_i\in N(v_1)$ and $y_j\in N(v_p)$ such that $y_j\prec x_i$ and $v_1v,v_pv\not\in E(G)$ for each vertex $v$ with $y_j\prec v\prec x_i$. Put
$$
C=v_1x_i\overrightarrow{P}v_py_i\overleftarrow{P}v_1.
$$
Then
\begin{equation*}
    \begin{split}
      c&\ge |V(C)|\ge |N(v_1)|+|N^+(v_p)|+|\{v_1\}|-|\{y_j\}| \\
&\ge|N(v_1)|+|N(v_p)|=d(v_1)+d(v_p). 
    \end{split}
\end{equation*}

Then we can argue as in Case 1.          \qquad \rule{7pt}{6pt}

\section{Concluding remarks}

The well-known Bondy-Chv\'atal’s power theorem covers all the previous theorems due to Dirac, Ore, Faudree, P\'osa, Bondy, Chv\'atal and Las-Vergnas. Often the same natural question arises: why do we need all these covered theorems, especially Dirac's theorem - the weakest and oldest among them, and why every time we go back and remember them? 

From a practical point of view, they are indeed outdated and unusable. But they are of valuable interest to scientific researchers, because each result that is at the origins or close to the origins as a starting point is a source for numerous generalizations and new ideas in a great  variety of ways. All these searches ultimately lead to results of practical value. 

In this regard, Theorem 1 and Theorem 2 proved in this paper, and all similar theorems are of particular interest.

\noindent Institute for Informatics and Automation Problems\\ 
National Academy of Sciences\\
P. Sevak 1, Yerevan 0014, Armenia\\
E-mail: zhora@iiap.sci.am

\end{document}